\documentclass[ 9pt]{amsart}
\usepackage{graphicx}
\usepackage{amssymb}
\usepackage{amsmath}
\usepackage{amsthm}
\usepackage{amscd}
\usepackage[all,2cell,dvips]{xy}
\UseAllTwocells \SilentMatrices
\newtheorem{thm}{Theorem}[section]

\newtheorem{cor}[thm]{Corollary}
\newtheorem{lem}[thm]{Lemma}

\newtheorem{prop}[thm]{Proposition}
\theoremstyle{definition}

\theoremstyle{remark}

\numberwithin{equation}{section}

\begin{document}
\title[Quotient Triangulated categories]{Quotient Triangulated categories}
\author[X. W. Chen,  \    P. Zhang
] {$^a$Xiao-Wu Chen, \ \  $^b$Pu Zhang$^*$}
\thanks{$^*$ The corresponding author}
\thanks{Supported by the National Natural Science Foundation of China (Grant No.
10301033) and of Shanghai City (Grant No. ZR0614049).}
%%\subjclass{}%
\thanks{E-mail: xwchen$\symbol{64}$mail.ustc.edu.cn,
pzhang$\symbol{64}$sjtu.edu.cn}
%\keywords{Quantum groups, quivers}%
\maketitle
\date{}%
\dedicatory{}%
\commby{}%
\begin{center}
$^A$Department of Mathematics, \ University of Science and
Technology of China \\Hefei 230026, P. R. China
\end{center}
\begin{center}
$^B$Department of Mathematics, \ \ Shanghai Jiao Tong University\\
Shanghai 200240, P. R. China\end{center}

\vskip 5pt

\begin{center}
Dedicated to Professor Zhe-Xian Wan on the occasion of his eightieth
birthday
\end{center}

\begin{abstract} For a self-orthogonal module $T$, the relation between
the quotient triangulated category $D^b(A)/K^b({\rm add} T)$ and the
stable category of the Frobenius category of $T$-Cohen-Macaulay
modules is investigated. In particular, for a Gorenstein algebra, we
get a relative version of the description of the singularity
category due to Happel. Also, the derived category of a Gorenstein
algebra is explicitly given, inside the stable category of the
graded module category of the corresponding trivial extension
algebra, via Happel's functor $F: D^b(A) \longrightarrow
T(A)^{\mathbb{Z}}\mbox{-}\underline{\rm
 mod}$.
\end{abstract}

\vskip20pt

\centerline {\bf Introduction}

\vskip10pt

Throughout, $A$ is a finite-dimensional associative algebra over a
field $k$, $A \mbox{-mod}$ the category of finite-dimensional left
$A$-modules,  and $D^b(A)$ the bounded derived category of
$A\mbox{-mod}$. Let $K^b(A\mbox{-inj})$ and $K^b(A\mbox{-proj})$ be
the bounded homotopy categories of injective, and projective
A-modules, respectively. Viewing them as thick (\'epaisse; see [V1])
triangulated subcategories of $D^b(A)$, one has the quotient
triangulated categories:
\begin{align*}
\mathcal{D}_I(A): = D^b(A)/{K^b(A\mbox{-inj})}\quad \mbox{and} \quad
\mathcal{D}_P(A):= D^b(A)/{K^b(A\mbox{-proj})}.
\end{align*}
Note that $\mathcal{D}_I(A) = 0$ (resp. $\mathcal{D}_P(A) = 0$) if
and only if ${\rm gl.dim A} < \infty$. In the same way, for an
algebraic variety $X$ one has the quotient triangulated category
${\bf D}_{Sg}(X): = {\bf D}^b(coh(X))/\textsf{perf}(X)$, where ${\bf
D}^b(coh(X))$ is the bounded derived category of coherent sheaves on
$X$, and $\textsf{perf}(X)$ is its full subcategory of perfect
complexes. Note that ${\bf D}_{Sg}(X) = 0$ if and only if $X$ is
smooth. Thus, for an algebra $A$ of infinite global dimension, or a
singular variety $X$,  it is of interest to investigate
$\mathcal{D}_I(A)$, $\mathcal{D}_P(A)$, and ${\bf D}_{Sg}(X)$. These
quotient triangulated categories become an important topic in
algebraic geometry and representation theory of algebras through the
work of  Buchweitz [Buc], Keller-Vossieck [KV1], Rickard [Ric2],
Happel [Hap2], Beligiannis [Bel], J${\o}$rgensen [J], Orlov [O1],
[O2], Krause [Kr], Krause-Iyengar [KI], and others, as they measure
the complexity of possible singularities. In particular, they are
called the singularity categories in [O1]; in [Ric2] (Theorem 2.1)
it was proved that for a self-injective algebra $A$,
$\mathcal{D}_P(A)$ is triangle-equivalent to the stable module
category of $A$ modulo the projectives; and in [Hap2] (Theorem 4.6)
it was proved that for a Gorenstein algebra $A$, $\mathcal{D}_P(A)$
is  triangle-equivalent to the stable category of the Frobenius
category of Cohen-Macaulay modules.

\vskip10pt

For a self-orthogonal $A$-module $T$ (i.e., ${\rm Ext}_A^i(T, T)=0$
for each $i\ge 1$), let ${\rm add} T$ denote the full subcategory of
$A\mbox{-mod}$ whose objects are the direct summands of finite
direct sum of copies of $T$. Then $K^b({\rm add} T)$ is a
triangulated subcategory of $D^b(A)$ ([Hap1], p.103). Since both
$K^b({\rm add} T)$ and $D^b(A)$ are Krull-Schmidt categories (see
[KV2], or [BD]), i.e., each object can be uniquely decomposed into a
direct sum of (finitely many) indecomposables and indecomposables
have local endomorphism rings, it follows that $K^b({\rm add} T)$ is
closed under direct summands, that is, $K^b({\rm add}T)$ is thick in
$D^b(A)$ (see Proposition 1.3 in [Ric2];  or [V2]), and hence one
has the quotient triangulated category $D^b(A)/K^b({\rm add} T).$ In
the view of the tilting theory (see e.g. [HR], [Rin1], [Hap1],
[AR1], [M]), $\mathcal{D}_I(A)$ (resp., $\mathcal{D}_P(A)$) is just
the special case of $\mathcal{D}_T(A)$ when $T$ is a generalized
cotilting module (resp., \ $T$ is a generalized tilting module).
This encourages us to look at $\mathcal{D}_I(A)$ and
$\mathcal{D}_P(A)$ in terms of generalized cotilting and tilting
modules, respectively, and to understand $\mathcal{D}_T(A)$ for
self-orthogonal modules $T$, in general. If $(A, {_AT_B}, B)$ is a
generalized (co)tilting triple, then by a theorem due to Happel
([Hap1], Theorem 2.10, p.109, for the finite global dimension case),
and due to Cline, Parshall, and Scott ([CPS], Theorem 2.1, for
general case. See also Rickard [Ric1], Theorem 6.4, in terms of
tilting complexes), this investigation permits us to understand the
singularity category $\mathcal{D}_P(B)$ in terms of ${_AT}$.

\vskip 10pt

For a self-orthogonal module $T$, we study in Section 2 the
relation between the quotient triangulated category
$D^b(A)/K^b({\rm add} T)$ and the stable category of the Frobenius
category $\alpha(T)$ of $T$-Cohen-Macauley modules (see 2.1 for
the definition of this terminology). In particular, for a
Gorenstein algebra, we get the relative version of the description
of the singularity category due to Happel ([Hap2], Theorem 4.6).
See Theorem 2.5.

\vskip 20pt

Denote by $T(A): = A\oplus D(A)$ the trivial extension of $A$, where
$D = {\rm Hom}_k(-, k)$. It is $\mathbb{Z}$-graded with ${\rm deg} A
= 0$ and ${\rm deg} D(A) = 1$. Denote by
$T(A)^{\mathbb{Z}}\mbox{-}{\rm
 mod}$ the category of finite-dimensional $\mathbb{Z}$-graded
 $T(A)$-modules with morphisms of degree $0$. This is a Frobenius
 abelian category, and hence its stable category $T(A)^{\mathbb{Z}}\mbox{-}\underline{\rm
 mod}$ modulo projectives is a triangulated category.
 Happel has constructed a fully faithful exact functor
 $F: D^b(A) \longrightarrow T(A)^{\mathbb{Z}}\mbox{-}\underline{\rm
 mod}$ ([Hap1], p.88, plus p.64);
 and $F$ is dense if and only ${\rm gl.dim A} <
 \infty$ ([Hap3]). Consider the natural embedding $ \emph{i}:
\mbox{A-mod}\hookrightarrow
 T(A)^{\mathbb{Z}}\mbox{-}\underline{\rm
 mod}$ such that each $A$-module $M$ is a graded $T(A)$-module
concentrated at degree $0$, which is the restriction of $F$. Denote
by $\mathcal{N}$, $\mathcal{M}_P$, and $\mathcal{M}_I$ the
triangulated subcategories of
$T(A)^{\mathbb{Z}}\mbox{-}\underline{\rm mod}$ generated by
$A\mbox{-}$mod, $A\mbox{-proj}$, and $A \mbox{-inj}$, respectively.
Since  both $T(A)^{\mathbb{Z}}\mbox{-}\underline{\rm mod}$ and
$\mathcal{N}$ are Krull-Schmidt, it follows that $\mathcal{N}$ is
closed under direct summands, that is, $\mathcal{N}$ is thick in
$T(A)^{\mathbb{Z}}\mbox{-}\underline{\rm mod}$; and so are
$\mathcal{M}_P$ and $\mathcal{M}_I$. So, one has the quotient
triangulated categories:
\begin{align*}
T(A)^{\mathbb{Z}}\mbox{-}\underline{\rm mod}/{\mathcal{N}}, \ \ \ \
T(A)^{\mathbb{Z}}\mbox{-}\underline{\rm mod}/\mathcal{M}_I, \ \ \ \
\mbox{and} \ \ \ T(A)^{\mathbb{Z}}\mbox{-}\underline{\rm
mod}/\mathcal{M}_P.
\end{align*}
Then Happel's theorem above reads as: there are equivalences of
triangulated categories
$$F: D^b(A) \simeq \mathcal{N}, \ \  \ \ F: K^b(A\mbox{-inj}) \simeq
\mathcal{M}_I, \ \ \ \ \mbox{and} \ \ \ F: K^b(A\mbox{-proj}) \simeq
\mathcal{M}_P;$$ and  $T(A)^{\mathbb{Z}}\mbox{-}\underline{\rm
mod}/{\mathcal{N}} = 0$ if and only if ${\rm gl.dim}\; A < \infty.$

\vskip10pt

It is then of interest to study these quotient triangulated
categories. In Section 3 we only take the first step by  giving an
explicit description of the bounded derived category of a Gorenstein
algebra $A$ inside $T(A)^{\mathbb{Z}}\mbox{-}\underline{\rm mod}$,
via Happel's functor above. See Theorem 3.1.

\vskip10pt

\section{\bf Preliminaries}

\subsection{} An algebra $A$ is Gorenstein if
${\rm proj.dim}\;  {_AD(A_A)}< \infty$ and ${\rm inj.dim} \; _AA <
\infty$. Self-injective algebras and algebras of finite global
dimension are Gorenstein; the tensor product $A\otimes_k B$ is
Gorenstein if and only if so are $A$ and $B$ ([AR2], Proposition
2.2). Note that $A$ is Gorenstein if and only if $K^b(A
\mbox{-proj})=K^b(A\mbox{-inj})$ inside $D^b(A)$ ([Hap2], Lemma
1.5). Thus, by Theorem 6.4 and Proposition 9.1 in [Ric1], if the
algebras $A$ and $B$ are derived equivalent, then $A$ is Gorenstein
if and only if so is $B$. Also, cluster-tilted algebras are
Gorenstein ([KR]).

\par \vskip 10pt

\subsection{} For basics on triangulated categories and derived
categories we refer to [Har] and [V1]. Following [BBD], the shift
functor in a triangulated category is denoted by $[1]$. Recall that,
by definition, triangulated subcategories are full subcategories. By
a multiplicative system, we will always mean a multiplicative system
compatible with the triangulation.  For a multiplicative system $S$
of a triangulated category $\mathcal K$, we refer to [Har] (see also
[V1] and [I]) for the construction of the quotient triangulated
category $S^{-1}\mathcal K$ via localization, in which morphisms are
given by right fractions (if one uses left fractions then one gets a
quotient triangulated category isomorphic to $S^{-1}\mathcal K$).

\vskip10pt

Let $\mathcal A$ be an abelian category, $\mathcal B$ a full
subcategory of $\mathcal A$,  and $\varphi: K^b(\mathcal
B)\longrightarrow D^b(\mathcal A)$ the composition of the embedding
$K^b(\mathcal B)\hookrightarrow K^b(\mathcal A)$ and the
localization functor $K^b(\mathcal A)\longrightarrow D^b(\mathcal
A)$. If $\varphi$ is fully faithful, then $K^b(\mathcal B)$ is a
triangulated subcategory of $D^b(\mathcal A)$.

Applying this to ${\rm add} T$ with $T$ a self-orthogonal
$A$-module, we know by Lemma 2.1 in [Hap1], p.103, that $K^b({\rm
add} T)$ is a triangulated subcategory of $D^b(A).$ If $T$ is a
generalized tilting module, then $K^b({\rm add} T) =
K^b(A\mbox{-proj})$ in $D^b(A)$ (in fact,  for any projective module
$P$ and $T'\in {\rm add}T$, we have $P\in K^b({\rm add} T)$ and
$T'\in K^b(A\mbox{-proj})$, in $D^b(A)$. Then the assertion follows
from the fact that $K^b(A\mbox{-proj})$ and $K^b({\rm add} T)$ are
the triangulated subcategories of $D^b(A)$ generated by ${\rm add}
A$ and by ${\rm add} T$, respectively).

\vskip10pt

\subsection {} An exact category $\mathcal{A}$ is a full subcategory of an abelian
category, closed under extensions and direct summands, together with
the exact structure given by the set of all the short exact
sequences of the ambient abelian category with terms in
$\mathcal{A}$. Such exact sequences will be referred as admissible
exact sequences (see Quillen [Q]; or conflations in the sense of
Gabriel-Roiter, see e.g. Appendix A in [K1]). A Frobenius category
is an exact category in which there are enough (relative) injective
objects and (relative) projective objects, such that the injective
objects coincide with the projective objects. For the reason
requiring that $\mathcal{A}$ is closed under direct summands see
Lemma 1.1 below. Compare p.10 in [Hap1], Appendix A in [K1], and
[Q]. Denote by $\underline{\mathcal{A}}$ its stable category. For a
morphism $u: \ X\longrightarrow Y$ in $\mathcal{A}$, denote its
image in $\underline{\mathcal{A}}$ by $\underline{u}$.

\vskip 10pt

\begin {lem} Let $\mathcal{A}$ be a
Frobenius category. Then $X\simeq Y$ in $\underline{\mathcal{A}}$ if
and only if there are injective objects $I$ and $J$ such that $X
\oplus J\simeq Y\oplus I$ in ${\mathcal{A}}$.\end{lem} \noindent
{\bf Proof.}\quad This is well-known. We include a proof for
convenience. Let $\underline f: X \longrightarrow Y$ be an
isomorphism in $\underline{\mathcal{A}}$. Then there exist an
injective object $I$,  $g: Y \longrightarrow X, $ $a: X
\longrightarrow I$ and $b: I \longrightarrow X$, such that $(g,
-b)\circ {f \choose a}={\rm Id}_X.$ Thus there exists $J$ in
$\mathcal{A}$ (here we need  $\mathcal{A}$ being closed under direct
summands) and $h: X\oplus J \simeq Y \oplus I$, such that ${f
\choose a} = h{1 \choose 0}$ in $\mathcal{A}$. Thus $\underline {{1
\choose 0}} : X \longrightarrow X\oplus J$ is an isomorphism in
$\underline {\mathcal{A}}$. Then by a matrix calculation we have
${\rm Id}_{\underline J} = 0$, i.e., $J$ is an injective object.
\hfill $\blacksquare$

\vskip10pt

Let $\mathcal{A}$ be a Frobenius category. Recall the triangulated
structure in $\underline{\mathcal{A}}$ from [Hap1], Chapter 1,
Section 2. The shift functor $[1]:
\underline{\mathcal{A}}\longrightarrow \underline{\mathcal{A}}$ is
defined such that for  $X\in \mathcal{A}$,
\begin{align}
0 \longrightarrow X \stackrel{i_X}{\longrightarrow} I(X)
\stackrel{\pi_X}{\longrightarrow} X[1] \longrightarrow 0
\end{align}
is an admissible exact sequence in $\mathcal{A}$ with $I(X)$ an
injective object. By Lemma 1.1, if $X\simeq Y$ in $\underline
{\mathcal{A}}$\ \ then  \ $X[1]\simeq Y[1]$ in
$\underline{\mathcal{A}}$, and the object $X[1]$ in $\underline
{\mathcal{A}}$ does not depend on the choice of $(1.1)$. For $u: \
X\longrightarrow Y$ in $\mathcal{A}$, the standard triangle $X
\stackrel{\underline{u}}{\longrightarrow} Y
\stackrel{\underline{v}}{\longrightarrow} C_u
\stackrel{\underline{w}}{\longrightarrow} X[1]$ \ in
$\underline{\mathcal{A}}$ \ is defined by the pushout diagram

\[
\xymatrix{ 0 \ar[r]  & X \ar[d]^u \ar[r]^-{i_X} & I(X) \ar[d]^{\bar
u} \ar[r]^-{\pi_X} & X[1]
\ar@{=}[d] \ar[r] & 0 \\
 0 \ar[r]  & Y \ar[r]^-{v} & C_u \ar[r]^-w & X[1] \ar[r] & 0;}
\]

\noindent and then the distinguished triangles in
$\underline{\mathcal{A}}$ \ are defined to be the triangles
isomorphic to the standard ones. We need the following fact in
[Hap1], p.22.

\vskip10pt

\begin{lem} Distinguished
triangles in $\underline{\mathcal{A}}$ are just given by short exact
sequences in $\mathcal{A}$.

More precisely, let $0 \longrightarrow X
\stackrel{u}{\longrightarrow} Y \stackrel{v}{\longrightarrow} Z
\longrightarrow 0$ be an admissible exact sequence in $\mathcal{A}$.
Then $X \stackrel{\underline{u}}{\longrightarrow} Y
\stackrel{\underline{v}}{\longrightarrow} Z
\stackrel{-\underline{w}}{\longrightarrow} X[1]$ is a distinguished
triangle in $\underline{\mathcal{A}}$, where $w$ is an $\mathcal
A$-map such that the following diagram is commutative (note that any
two such maps $w$ and $w'$ give the isomorphic triangles)
\begin{align}
\xymatrix{ 0 \ar[r]  & X \ar@{=}[d] \ar[r]^-u & Y \ar@{.>}[d]^\sigma \ar[r]^-v & Z \ar@{.>}[d]^-w \ar[r] & 0 \\
 0 \ar[r]  & X \ar[r]^-{i_X} & I(X) \ar[r]^-{\pi_X} & X[1] \ar[r] & 0.}
\end{align}

\vskip5pt

Conversely, let  $X' \stackrel{\underline{u'}}{\longrightarrow} Y'
\stackrel{\underline{v'}}{\longrightarrow} Z'
\stackrel{-\underline{w'}}{\longrightarrow} X'[1]$ be a
distinguished triangle in $\underline{\mathcal{A}}$.  Then there is
an admissible exact sequence $0 \longrightarrow X
\stackrel{u}{\longrightarrow} Y \stackrel{v}{\longrightarrow} Z
\longrightarrow 0$ in $\mathcal{A}$, such that the induced
distinguished triangle  $X \stackrel{\underline{u}}{\longrightarrow}
Y \stackrel{\underline{v}}{\longrightarrow} Z
\stackrel{-\underline{w}}{\longrightarrow} X[1]$ is isomorphic to
the given one, where $w$ is an $\mathcal{A}$-map such that
\emph{(1.2)} is commutative.\end{lem}

\noindent {\bf Proof.} \ We include a proof of the second part for
convenience. Let $X' \stackrel{\underline{u'}}{\longrightarrow} Y'
\stackrel{\underline{v'}}{\longrightarrow} Z'
\stackrel{-\underline{w'}}{\longrightarrow} X'[1]$ be a
distinguished triangle in $\underline{\mathcal{A}}$. Then it is
isomorphic to a standard triangle $X
\stackrel{\underline{u}}{\longrightarrow} Y
\stackrel{\underline{v}}{\longrightarrow} C_u
\stackrel{\underline{w}}{\longrightarrow} X[1]$  in
$\underline{\mathcal{A}}$, with an admissible exact sequence in
${\mathcal{A}}$:
\begin{align}0\longrightarrow
X \stackrel{\binom u {i_X}} {\longrightarrow} Y\oplus I(X)
\stackrel{(v, \ -\bar u)}{\longrightarrow} C_u \longrightarrow
0\end{align} and the  commutative diagram

\[
\xymatrix{ 0 \ar[r]  & X \ar@{=}[d] \ar[r]^-{\binom u {i_X}} &
Y\oplus I(X) \ar[d]^{(0,1)} \ar[r]^-{(v, \; -\bar u)} & C_u
\ar[d]^-{-w} \ar[r]
& 0 \\
 0 \ar[r]  & X \ar[r]^-{i_X} & I(X) \ar[r]^-{\pi_X} & X[1] \ar[r] & 0.}
\]
This implies that the distinguished triangle induced by (1.3) is
isomorphic to the standard triangle. \hfill $\blacksquare$

\vskip10pt

\begin{lem}  Let $\mathcal{A}$ be a  Frobenius category.
Then there is a bijection between the class of the full
subcategories $\mathcal B$ of ${\mathcal{A}}$, where $\mathcal B$
contains all the injective objects of $\mathcal{A}$, such that if
two terms in an admissile exact sequence in $\mathcal A$ lie in
$\mathcal B$, then the third term also lies in $\mathcal B$, and the
class of triangulated subcategories of $\underline{\mathcal{A}}$.
\end{lem}
\noindent {\bf Proof.}\quad If $\mathcal{B}$ is such a full
subcategory of $\mathcal{A}$, then by Lemmas 1.2 and 1.1
$\underline{\mathcal{B}}$ is a triangulated subcategory of
$\underline{\mathcal{A}}$. Conversely, let $\mathcal D$ be a
triangulated subcategory of $\underline{\mathcal{A}}$. Set
$$\mathcal B: = \{X\in \mathcal{A}\ | \  \mbox{there exists} \ Y\in \mathcal
D \ \mbox{such that}\ X\simeq Y \ \mbox{in} \ \underline{\mathcal A}
\}.$$ Then $\mathcal D = \underline {\mathcal B}$. Since $\mathcal
D$ contains zero object, it follows that $\mathcal B$ contains all
the injective objects of $\mathcal A$; and by Lemma 1.1 $\mathcal B$
has the required property. If $\mathcal{B}$ and $\mathcal{B'}$ are
two such a different full subcategories  of $\mathcal{A}$, then by
Lemma 1.1 $\underline{\mathcal{B}}$ and $\underline{\mathcal{B'}}$
are also different. \hfill $\blacksquare$

\vskip15pt

\section{\bf Quotient triangulated category $D^b(A)/K^b({\rm add} T)$}

\vskip10pt

Throughout the section, $T$ is a self-orthogonal $A$-module. All
subcategories will be assumed to be closed under isomorphisms and
finite direct sums. However, following [Rin2], we do not assume that
they are closed under direct summands. \vskip10pt

\subsection{} Consider the following full subcategories of
$A\mbox{-mod}$ introduced in [AR1]:
\begin{align*}
T^\perp:=\{X \ | \ {\rm Ext}_A^i(T, X)=0, \ \forall \ i\ge 1\},
\end{align*}
\begin{align*}
{\rm add}^\sim {T}:= \{X \ |\  \exists \ \mbox {an exact
sequence}\ \cdots \longrightarrow T^{-i} \longrightarrow \cdots
  T^0 \longrightarrow X \longrightarrow 0,  \ T^{-i}\in {\rm add}T,  \forall \
  i\},
\end{align*}
\begin{align*}
_T\mathcal X : =  \{X \ |\  \exists \ \mbox {an exact sequence}\ \ &
\cdots \longrightarrow T^{-i}\stackrel{d^{-i}}{\longrightarrow}
T^{-(i-1)} \longrightarrow \cdots  \stackrel{d^{-1}}
{\longrightarrow}
  T^0 \stackrel{d^0} \longrightarrow X {\longrightarrow} 0,  \\ &
  T^{-i}\in {\rm add}T, \ \    {\rm Ker} d^{-i}\in \ T^\perp,\ \forall \
  i\ge 0\},
\end{align*}
\noindent and
\begin{align*} {\rm add} ^ {\wedge} T:= \{X \ | \
\exists \ \mbox {an exact sequence}\  0 \longrightarrow T^{-n}
\longrightarrow \cdots
  T^0 \longrightarrow X \longrightarrow 0,  \ T^{-i}\in {\rm add}T, \   \ \forall \
   i \ \}.
\end{align*}

By dimension-shifting we have \ ${\rm add}^{\wedge} T \subseteq \
_T\mathcal X \subseteq {\rm add}^\sim T\ \cap \ T^\perp.$ Note that
$\ _T\mathcal X = {\rm add}^\sim T$\  if and only if \ ${\rm
add}^{\sim} T \subseteq {T^{\perp}}$ (For the ``if part", note that
${\rm Ker} d^{-i}$ is still in ${\rm add}^\sim T$, and hence in
$T^\perp$). If $T$ is exceptional (i.e., ${\rm proj.dim}\; T<
\infty$ and $T$ is self-orthogonal), then ${\rm add}^\sim T
\subseteq T^{\perp}$, and hence $\ _T\mathcal X = {\rm add}^\sim T$.

\par \vskip 10pt

If $T$ is a generalized tilting module (i.e., $T$ is exceptional,
and there is an exact sequence $0\longrightarrow
{_AA}\longrightarrow T^0\longrightarrow T^1 \longrightarrow \cdots
\longrightarrow T^n\longrightarrow 0$ with each $T^i\in {\rm add}
T$), then ${\rm add}^\sim T = T^{\perp}$, and hence $\ _T\mathcal X
= {\rm add}^\sim T = T^{\perp}$. (In fact, by the theory of
generalized tilting modules, $X\in T^{\perp}$ can be generated by
$T$, see [M], Lemma 1.8; and then by using a classical argument in
[HR], p.408, one can prove $X\in {\rm add}^\sim T$ by induction.)

\par \vskip 10pt

If ${\rm gl.dim}A < \infty$ then ${\rm add}^\wedge T = \ _T\mathcal
X = {\rm add}^\sim T$ for any self-orthogonal module $T$. For this
it suffices to prove $_T\mathcal X\subseteq {\rm add}^\wedge T$.
This follows from
$$\operatorname{Ext}^1_A({\rm Ker} d^{-(i-1)}, {\rm Ker} d^{-i}) \cong
\operatorname{Ext}^2_A({\rm Ker} d^{-(i-2)}, {\rm Ker} d^{-i})
\cong\cdots \cong \operatorname{Ext}^i_A({\rm Ker} d^{0}, {\rm Ker}
d^{-i}) = 0,$$ where $i\gg 0$.

\vskip10pt

Dually, we have the concept of a generalized cotilting module,
full subcategories of $A\mbox{-mod}$: \ $^{\perp}T, \ {\rm
add}_\sim T\ , \ \mathcal X_T, \ {\rm add}_\wedge T$, and the
corresponding facts.

\vskip10pt

We define the category of $T$-Cohen-Macaulay modules to be the
subcategory $\mathfrak{a}(T): = \mathcal{X}_T \cap { _T
\mathcal{X}}$. One may compare it with the definition of the
category Cohen-Macaulay modules given in [Hap2] and [AR2]. By
Proposition 5.1 in [AR1], both subcategories $\mathcal{X}_T$ and $_T
\mathcal{X}$ are closed under extensions and direct summands, and
thus so is $\alpha(T)$. One can easily verify that $\mathfrak{a}(T)$
is a Frobenius category, where ${\rm add}T$ is exactly the full
subcategory of all the (relatively) projective and injective
objects.

\vskip10pt

\subsection {} For $M, N\in A\mbox{-mod}$, let $T(M,
N)$ denote the subspace of $A$-maps from $M$ to $N$ which factor
through ${\rm add}T$. The following lemma seems to be of independent
interest. We give an explicit proof by using calculus of fractions.

\vskip10pt

\begin{lem} Let $T$ be a self-orthogonal module. If $M\in \mathcal X_T$ and $N\in T^\perp$, then there is a
canonical isomorphism of $k$-spaces $${\rm Hom}_A(M, N)/{T(M,
N)}\simeq {\rm Hom}_{D^b(A)/{K^b({\rm add} T)}} (M, N).$$
\end{lem}
\noindent {\bf Proof}.\ In what follows, a doubled arrow means a
morphism belonging to the saturated multiplicative system,
determined by the thick triangulated subcategory ${K^b({\rm add}T)}$
of $D^b(A)$ (see [V1], [Har], or [I]). A morphism from $M$ to $N$ in
$D^b(A)/{K^b({\rm add}T)}$ is denoted by right fraction $a/s: \ M
\stackrel{s}{\Longleftarrow} Z^\bullet
\stackrel{a}{\longrightarrow}N$, where $Z^\bullet \in D^b(A)$. Note
that the mapping cone ${\rm Con}(s)$ lies in $K^b({\rm add}T)$. We
have a distinguished triangle in $D^b(A)$

\begin{align}
 Z^\bullet \stackrel{s}{\Longrightarrow} M
\longrightarrow {\rm Con}(s) \longrightarrow Z^\bullet[1].
\end{align}

\vskip5pt

Consider the $k$-map $G: {\rm Hom}_A(M, N) \longrightarrow {\rm
Hom}_{D^b(A)/{K^b({\rm add}T)}}(M, N)$,  given by $G(f)=f/{{\rm
Id}_M}$. First, we prove that $G$ is surjective. By $M\in \mathcal
X_T$ we have an exact sequence
\begin{align*}
0\longrightarrow M \stackrel{\varepsilon}{\longrightarrow}
T^0\stackrel{d^0}{\longrightarrow} T^1
\stackrel{d^1}{\longrightarrow} \cdots \longrightarrow T^{n}
\stackrel{d^n}{\longrightarrow} \cdots
\end{align*}
with ${\rm Im} d^{i}\in \ ^\perp T, \ \forall \ i\ge 0.$ Then $M$
is isomorphic in $D^b(A)$ to the complex $T^\bullet: =
0\longrightarrow T^0 \longrightarrow T^1 \longrightarrow \cdots,$
and then isomorphic to the complex $0\longrightarrow T^0
 \longrightarrow \cdots \longrightarrow T^{l-1}
\longrightarrow {\rm Ker}d^{l}\longrightarrow 0$ for each $l\ge
1$. The last complex induces a distinguished triangle in $D^b(A)$

\begin{align}\sigma^{<
l} T^\bullet [-1] \longrightarrow {\rm
Ker}d^{l}[-l]\stackrel{s'}{\Longrightarrow} M
\stackrel{\varepsilon}{\longrightarrow} \sigma^{< l}
T^\bullet,\end{align}

\vskip5pt

\noindent where $\sigma^{<l}T^\bullet=0 \longrightarrow T^0
\longrightarrow T^1 \longrightarrow \cdots \longrightarrow T^{l-1}
\longrightarrow 0$, and the mapping cone of $s'$ lies in $K^b({\rm
add}T)$. Since ${\rm Ker}d^l \in {^\perp T}$ and  ${\rm Con}(s)\in
K^b({\rm add}T)$, it follows that there exists $l_0\gg 0$ such that
for each $l\geq l_0$
\begin{align*}
{\rm Hom}_{D^b(A)}({\rm Ker}d^l[-l], {\rm Con}(s))=0.
\end{align*}

(To see this, let ${\rm Con}(s)$ be of the form $0\longrightarrow
W^{-t'}\longrightarrow \cdots \longrightarrow W^t\longrightarrow 0$
with $t', t\ge 0,$ and each $W^i\in {\rm add}T$. Consider the
distinguished triangle in $D^b(A)$

$$\sigma^{< t} {\rm Con}(s)[-1]\longrightarrow W^t[-t]\longrightarrow {\rm Con}(s)
\longrightarrow \sigma^{ < t }{\rm Con}(s),$$

\vskip5pt

\noindent
 Take $l_0$ to be $t+1$, and apply the functor ${\rm
Hom}_{D^b(A)} ({\rm Ker}d^{l}[-l], -)$ to this distinguished
triangle. Then the assertion follows from  ${\rm Ker}d^l \in
{^\perp T}$ and induction.) \par

\vskip5pt

Write $E={\rm Ker}d^{l_0}$, and take $l=l_0$ in (2.2). By applying
${\rm Hom}_{D^b(A)}(E[-l_0], -)$ to (2.1) we get $h: E[-l_0]
\longrightarrow Z^\bullet$ such that $s'=s \circ h$. So we have
$a/s=(a\circ h)/s'$.\par

Apply ${\rm Hom}_{D^b(A)}(-, N)$ to (2.2), we get an exact sequence

$${\rm Hom}_{D^b(A)}(M, N) \longrightarrow {\rm
Hom}_{D^b(A)}(E[-l_0], N) \longrightarrow {\rm
Hom}_{D^b(A)}(\sigma^{<l_0}T^\bullet[-1], N).$$

\noindent We claim that ${\rm
Hom}_{D^b(A)}(\sigma^{<l_0}T^\bullet[-1], N)={\rm Hom}_{D^b(A)}
(\sigma^{<l_0}T^\bullet, N[1])=0$.

(In fact, apply ${\rm Hom}_{D^b(A)}(-, N[1])$ to the following
distinguished triangle in $D^b(A)$
\begin{align*}
\sigma^{<l_0-1}T^\bullet [-1]\longrightarrow T^{l_0-1}[1-l_0]
\longrightarrow \sigma^{<l_0}T^\bullet \longrightarrow
\sigma^{<l_0-1} T^\bullet.
\end{align*}
Then the assertion follows from induction and the assumption $N
\in T^\perp$.)

Thus, there exists $f: M \longrightarrow N$ such that $f \circ s'=a
\circ h$. So we have $a/s=(a\circ h)/s'=(f\circ s')/s'=f/{{\rm
Id}_M}$. This shows that $G$ is surjective.\par

\vskip5pt

On the other hand, if $f: M \longrightarrow N$ with $G(f)=f/{{\rm
Id}_M}=0$ in $D^b(A)/{K^b({\rm add}T)}$, then there exists $s:
Z^\bullet \Longrightarrow M$ with ${\rm Con}(s) \in K^b({\rm add}T)$
such that $f\circ s=0$. Use the same notation as in (2.1) and (2.2).
By the  argument above we have $s'=s \circ h$, and hence $f\circ
s'=0$. Therefore, by applying ${\rm Hom}_{D^b(A)} (-, N)$ to (2.2)
we see that there exists $f': \sigma^{<l_0}T^\bullet \longrightarrow
N$ such that $f' \circ \varepsilon=f$.\par

Consider the following distinguished triangle in $D^b(A)$
\begin{align*}
T^0[-1]\longrightarrow \sigma^{>0} (\sigma^{<l_0}) T^\bullet
\longrightarrow \sigma^{<l_0} T^\bullet
\stackrel{\pi}{\longrightarrow} T^0,
\end{align*}
where $\sigma^{>0} (\sigma^{<l_0}) T^\bullet =0 \longrightarrow
T^1 \longrightarrow T^2 \longrightarrow \cdots \longrightarrow
T^{l_0-1} \longrightarrow 0$, and $\pi$ is the natural morphism.
Again since $N\in T^\perp $, it follows that ${\rm Hom}_{D^b(A)}
(\sigma^{>0} (\sigma^{<l_0}) T^\bullet, N)=0$. By applying ${\rm
Hom}_{D^b(A)} (-, N)$ to the above triangle we obtain an exact
sequence
\begin{align*}
{\rm Hom}_{D^b(A)} (T^0, N) \longrightarrow {\rm Hom}_{D^b(A)} (
\sigma^{<l_0} T^\bullet , N) \longrightarrow 0.
\end{align*}
It follows that there exists $g: T^0 \longrightarrow N$ such that
$g\circ \pi=f'$. Hence $f=g\circ (\pi\circ \varepsilon)$. Since
$A\mbox{-mod}$ is a full subcategory of $D^b(A)$, it follows that
$f$ factors through $T^0$ in $A\mbox{-mod}$. This proves that the
kernel of $G$ is $T(M, N)$, which completes the proof. \hfill
$\blacksquare$

\par

\vskip10pt

Consider the natural functor $\mathcal X_T \cap T^{\perp}
\longrightarrow D^b(A)/{K^b({\rm add}T)}$, which is the composition
of the  embedding  $\mathcal X_T \cap T^{\perp} \hookrightarrow
A\mbox{-mod}$ and the  embedding $A\mbox{-mod} \hookrightarrow
D^b(A)$, and the localization functor $D^b(A) \longrightarrow
D^b(A)/{K^b({\rm add}T)}.$ Let $\underline {^{\perp}T \cap
T^{\perp}}$ denote the stable category of $^{\perp}T\cap T^{\perp}$
modulo ${\rm add} T$.

\vskip10pt

\begin{lem} Let $T$ be a generalized cotilting
module. Then the natural functor ${^\perp T} \ \cap \
T^{\perp}\longrightarrow D^b(A)/{K^b({\rm add}T)}$  \ induces a
fully faithful functor $$\underline {{^\perp T}\ \cap \
T^{\perp}}\longrightarrow D^b(A)/{K^b({\rm add} T)} =
\mathcal{D}_I(A).$$
\end{lem}
\noindent {\bf Proof}.\ Since $T$ is generalized cotilting then
$\mathcal X_T = {^{\perp} T}$, it follows that $K^b({\rm add} T) =
K^b(A \mbox{-inj})$ in $D^b(A)$. So the assertion follows from Lemma
2.1. \hfill $\blacksquare$

\vskip10pt

\subsection{} Let $T$ be a self-orthogonal module. By Lemma 2.1 the natural functor
$\underline{\mathcal X_T \cap T^{\perp}} \longrightarrow
D^b(A)/{K^b({\rm add}T)}$ is fully faithful. It is of interest to
know when it is dense. \vskip10pt

\begin{lem} $(i)$ Assume that ${\rm inj.dim}\;
_AA<\infty$. Let $T$ be a generalized cotilting $A$-module. Then the
natural functor $^\perp T \longrightarrow D^b(A)/{K^b({\rm add}T)} =
\mathcal{D}_I(A)$ is dense.\par

\vskip5pt

$(ii)$ If $A$ is Gorenstein and $T$ is a generalized tilting module,
then the natural functor \
$$^\perp T \cap T^{\perp} \longrightarrow
D^b(A)/{K^b({\rm add}T)} = \mathcal{D}_I(A)= \mathcal{D}_P(A)$$ is
dense.
\end{lem}
\noindent {\bf Proof}.\quad $(i)$\ \ By Happel (the dual of Lemma
4.3 in [Hap2]), the natural functor $A\mbox{-mod}\longrightarrow
\mathcal{D}_I(A)$ is dense. So for any object $X^\bullet\in
\mathcal{D}_I(A)$, there exists a module $M$ such that $X^\bullet
\cong M$ in $\mathcal{D}_I(A)$. Consider the subcategory
$$\widehat{\mathcal X_T}:=\{M \; |\; \exists \ \mbox
{an exact sequence}\ 0 \longrightarrow X^{-n} \longrightarrow \cdots
  X^0 \longrightarrow M \longrightarrow 0,  \ X^{-i}\in \mathcal X_T,  \forall \
  i \}.$$ Since $T$ is
generalized cotilting, it follows that $\widehat{\mathcal X_T} =
A\mbox{-mod}$, and hence by Auslander and Buchweitz (Theorem 1.1 in
[AB]), there exists an exact sequence
\begin{align}0\longrightarrow Y_M\longrightarrow X_M \longrightarrow
M \longrightarrow 0\end{align} with $Y_M\in {\rm add}^\wedge T, \
X_M\in \mathcal X_T = \ ^\perp T$. This induces a triangle
$Y_M\longrightarrow X_M \longrightarrow M \longrightarrow Y_M[1]$
in $D^b(A)$, and hence a triangle in $\mathcal{D}_I(A)$. Since $T$
is of finite injective dimension, it follows that modules in ${\rm
add}^\wedge T$ are of finite injective dimensions, and hence $Y_M
= 0$ in $\mathcal{D}_I(A)$, which implies $M\cong X_M$ in
$\mathcal{D}_I(A)$. This proves $(i)$.

\vskip5pt

$(ii)$ \ \ Since $A$ is Gorenstein it follows from Lemma 1.5 in
[Hap2] that $D^b(A)/{K^b({\rm add}T)} = \mathcal{D}_I(A)=
\mathcal{D}_P(A)$. For any object $X^\bullet\in \mathcal{D}_P(A)$,
by the dual of $(i)$, there exists a module $M\in T^\perp$ such
that $X^\bullet \cong M$ in $\mathcal{D}_P(A)$. Since $A$ is
Gorenstein, it follows from Lemma 1.3 in [HU] that $T$ is also
generalized cotilting. Repeat the argument in the proof of $(i)$
and note that ${\rm add}^\wedge T\subseteq T^\perp$.  It follows
from $(2.3)$ that $X_M\in T^\perp$, and hence $X^\bullet \cong
M\cong X_M\in \ ^\perp T \cap T^{\perp}$. This completes the
proof. \hfill $\blacksquare$

\vskip10pt

\subsection{} Let $T$ be a self-orthogonal $A$-module. Recall from 2.1  the category $\alpha(T)$
of $T$-Cohen-Macaulay modules, which is a Frobenius category and
${\rm add}T$ is exactly the full subcategory of all the (relatively)
projective and injective objects. Then the stable category of
$\mathfrak{a}(T)$ modulo ${\rm add}T$, denoted by
$\underline{\mathfrak{a}(T)}$, is a triangulated category. Since
$\mathfrak{a}(T)$ is a full subcategory of $\mathcal{X}_T \cap
T^\perp$, it follows from Lemma 2.1 that we have a natural fully
faithful functor $\underline{\mathfrak{a}(T)} \longrightarrow
D^b(A)/{K^b({\rm add}T)}$.

\vskip10pt

\begin{lem}\ \ Let  $T$ be a self-orthogonal module. Then
the natural embedding $\underline{\mathfrak{a}(T)} \longrightarrow
D^b(A)/{K^b({\rm add}T)}$ is an exact functor.
\end{lem}
\noindent {\bf Proof}. \ This follows from 1.2 in [K2]. We include a
direct justification.

Let $0 \longrightarrow X \stackrel{u}{\longrightarrow} Y
\stackrel{v}{\longrightarrow} Z \longrightarrow 0$ be an exact
sequence in $\mathfrak{a}(T)$, and $0\longrightarrow X
\stackrel{i_X}{\longrightarrow} T(X)
\stackrel{\pi_X}{\longrightarrow} S(X) \longrightarrow 0$ an exact
sequence with $T(X)\in {\rm add}T$ and $S(X)\in \mathfrak{a}(T)$.
Then $X \stackrel{\underline {u}}{\longrightarrow} Y
\stackrel{\underline {v}}{\longrightarrow} Z \stackrel{-\underline
{w}}{\longrightarrow} S(X)$ is a distinguished triangle in
$\underline{\mathfrak{a}(T)}$, where $w$ is an $A$-map such that the
following diagram is commutative
\begin{align}
\xymatrix{
0 \ar[r] & X \ar[r]^-u \ar@{=}[d] & Y\ar[d]^-\rho \ar[r]^-{v} \ar[r] & Z \ar[d]^-w \ar[r] &0 \\
0 \ar[r]  & X        \ar[r]^-{i_X} & T(X) \ar[r]^-{\pi_X}  & S(X)
\ar[r] & 0;
 }\end{align}
Note that any distinguished triangle in
$\underline{\mathfrak{a}(T)}$ is given in this way (cf. Lemma 1.2).

On the other hand, we have a distinguished triangle in $D^b(A)$
\begin{align}X \stackrel{u}{\longrightarrow} Y
\stackrel{v}{\longrightarrow} Z \stackrel{w'}{\longrightarrow}
X[1]\end{align} with \begin{align}w'= p_X/{v'}\in {\rm
Hom}_{D^b(A)}(Z, X[1])\end{align} as right fractions, where $p_X:
{\rm Con}(u) \longrightarrow X[1]$ is the natural morphism of
complexes,  and $v': {\rm Con}(u) \longrightarrow Z$ is the
quasi-isomorphism induced by $v$. Denote by  $p'_X: {\rm Con}(i_X)
\longrightarrow X[1]$ the natural morphism of complexes,  and
$\pi'_X: {\rm Con}(i_X) \longrightarrow S(X)$ the quasi-isomorphism
induced by $\pi_X$. Then right fraction
 $\beta_X:= - p'_X/{\pi'_X}$ is in ${\rm Hom}_{D^b(A)}(S(X), X[1])$.
We claim that $w'= - \beta_X w$ in $D^b(A)$, and hence  by $(2.5)$,
$X \stackrel{u} \longrightarrow Y \stackrel{v} \longrightarrow Z
\stackrel{-\beta_X w}\longrightarrow X[1]$ is a distinguished
triangle in $D^b(A)$, and hence  a distinguished triangle in
 $D^b(A)/{K^b({\rm add}T)}$.\par

In fact,  by $(2.6)$ the claim is equivalent to $p_X = - \beta_X (w
v')$ in $D^b(A)$. Denote by $\rho'$ the chain map ${\rm Con}(u)
\longrightarrow {\rm Con}(i_X)$  induced by $\rho$. Then $- \beta_X
(w v') = (p'_X/\pi'_X)(wv') = p'_X\rho' = p_X,$ where the second
equality follows from the multiplication rule of right fractions and
$wv=\pi_X\rho$ in $(2.4)$.

By the distinguished triangle $X \stackrel{i_X}{\longrightarrow}
T(X) \longrightarrow {\rm Con}(i_X)\stackrel{p'_X} \longrightarrow
X[1]$ in $D^b(A)$ we get the corresponding one in $D^b(A)/{K^b({\rm
add}T)}$ with $T(X) = 0$, it follows that $p'_X$, and hence
$\beta_X$, is an isomorphism in $D^b(A)/{K^b({\rm add}T)}$. This
shows that $\beta: G\circ S\longrightarrow [1]\circ G$ is a natural
isomorphism, where $G$ denote the natural functor
$\underline{\mathfrak{a}(T)} \longrightarrow D^b(A)/{K^b({\rm
add}T)}$. This completes the proof. \hfill $\blacksquare$

\vskip10pt

\subsection {} The following result gives a relative version of the explicit description due to
Happel, of the singularity categories of Gorenstein algebras
(Theorem 4.6 in [Hap2]).

\vskip10pt

\begin{thm} \ Let $A$ be a Gorenstein algebra and $T$ a generalized tilting $A$-module.
Then the natural functor induces a triangle-equivalence $\underline
{^{\perp}T \ \cap T^{\perp}}\cong D^b(A)/{K^b({\rm add}T)} =
\mathcal{D}_P(A) = \mathcal{D}_I(A).$
\end{thm}
\noindent {\bf Proof}.\quad  This follows from Lemmas 2.2, 2.3$(ii)$
and 2.4. Note that for a Gorenstein algebra, generalized tilting
modules coincide with generalized cotilting modules (Lemma 1.3 in
[HU]); and in this case $^{\perp}T \ \cap T^{\perp} = \mathcal{X}_T
\cap { _T \mathcal{X}}$. \hfill $\blacksquare$

\vskip10pt

Let us remark that one can also use the derived equivalence given by
$T$ to reduce Theorem 2.5 to the classical one where $T=A$ as in
Theorem 4.6 of [Hap2].

\vskip10pt

\subsection {} The following result is different from Lemma 2.3, and
seems to be of interest.

\vskip10pt

\begin{prop} Assume that ${\rm inj.dim}\;
_AA<\infty$. Let $T$ be a generalized tilting module. Then the
natural functor $^\perp T \cap T^\perp \longrightarrow
\mathcal{D}_I(A)$ is dense.\par
\end{prop}

\noindent {\bf Proof}.  \ \ Set $t: = {\rm proj.dim}\; T$. Since $T$
is generalized tilting and ${\rm inj.dim}\; _AA <\infty$, it follows
that $K^b({\rm add}T) = K^b(A\mbox{-proj}) \subseteq
K^b(A\mbox{-inj}),$ and hence ${\rm inj.dim}\; T = s < \infty$.

\vskip10pt

Identify $D^b(A)$ with $K^{+, b}(A\mbox{-inj})$. For any object
$I^\bullet$ in $ \mathcal{D}_I(A)$,  without loss of generality, we
may assume that
$$I^\bullet = 0 \longrightarrow I^0 \longrightarrow \cdots \longrightarrow
I^{l-1}  \longrightarrow I^l \stackrel {d^l}\longrightarrow \cdots
\longrightarrow I^{l+r-1}\longrightarrow I^{l+r}  \stackrel
{d^{l+r}} \longrightarrow \cdots$$ with $H^n(I^\bullet) = 0$ for
$n\ge l.$ Set $E: = {\rm Ker} d^{l+r}$ and $X: = {\rm Ker} d^{l}$.
Then the complex $0 \longrightarrow I^0 \longrightarrow \cdots
\longrightarrow I^{l-1} \longrightarrow I^l \stackrel
{d^l}\longrightarrow \cdots \longrightarrow I^{l+r-1}\longrightarrow
E\longrightarrow 0$ is quasi-isomorphic to $I^\bullet$, and hence
$I^\bullet \simeq E[-(l+r)]$ in $\mathcal D_I(A)$.

\vskip10pt

Take $r\ge s, t$. By the exact sequence of $A$-modules

$$0\longrightarrow X \longrightarrow I^l \longrightarrow \cdots
\longrightarrow I^{l+r-1} \longrightarrow E\longrightarrow 0,$$ and
${\rm proj.dim}\; T  = t <\infty$  and $r\ge t$, we infer that $E\in
T^\perp$.

By the generalized tilting theory we have $T^\perp = {\rm
add}^\sim T = \ _T\mathcal X$, and hence we have an exact sequence
of $A$-module

\begin{align}
0\longrightarrow W \longrightarrow T^{-(l+r-1)} \longrightarrow
\cdots \longrightarrow T^0 \longrightarrow E\longrightarrow 0
\end{align}
with each $T^i\in {\rm add}T$ and $W\in T^\perp$. Since  $K^b({\rm
add}T) \subseteq K^b(A\mbox{-inj}),$ it follows that the complex $0
\longrightarrow T^{-(l+r-1)} \longrightarrow \cdots \longrightarrow
T^0 \longrightarrow 0$ is in $K^b(A\mbox{-inj}),$ and hence $E =
W[l+r]$ in $\mathcal{D}_I(A)$. Since $T$ is self-orthogonal with
${\rm inj.dim}\; T = s$ and $r\ge s$, by $(2.7)$ we infer that $W\in
\ ^\perp T$. Thus $I^\bullet = W$ in $\mathcal{D}_I(A)$ with $W \in
\ ^\perp T\cap T^\perp$. \quad \hfill $\blacksquare$

\vskip10pt

\begin{cor} \ \ The following are equivalent

\vskip5pt

$(i)$ \ \ $\rm {gl.dim} A < \infty$;

\vskip5pt

$(ii)$ \ \ ${\rm inj.dim}\; _AA < \infty$, and \ $^{\perp}T \cap \
T^{\perp} = {\rm add} T$ for any generalized tilting module.
\end{cor}

\noindent {\bf Proof}. The implication of $(i) \Longrightarrow (ii)$
follows from Theorem 2.5; and  $(ii) \Longrightarrow (i)$ follows
from Proposition 2.6, since $K^b({\rm add}T) = K^b(A\mbox{-proj})
\subseteq K^b(A\mbox{-inj})$. \quad \hfill $\blacksquare$

\vskip 15pt

\section{\bf Bounded derived categories of Gorenstein algebras}

\vskip 10pt

\subsection{}Keep the notation in the introduction throughout this section, in particular for
$\mathcal{N}$ and $\mathcal{M}_P$.  An object in
$T(A)^{\mathbb{Z}}\mbox{-mod}$ and in
$T(A)^{\mathbb{Z}}\mbox{-}\underline{\rm mod}$ is denoted by
$M=\oplus_{n \in \mathbb{Z}} M_{n}$ with each $M_{n}$ an $A$-module
and $D(A).M_{n} \subseteq M_{n+1}$. The following result explicitly
describes the bounded derived category of a Gorenstein algebra $A$
inside $T(A)^{\mathbb{Z}}\mbox{-}\underline{\rm mod}$. \vskip10pt

\begin{thm}  Let $A$ be a Gorenstein algebra. Then via
Happel's functor $F: D^b(A) \longrightarrow
T(A)^{\mathbb{Z}}\mbox{-}\underline{\rm
 mod}$ we have
$$D^b(A)\simeq \mathcal{N}=\{\oplus_{n \in \mathbb{Z}} M_{n}\in T(A)^{\mathbb{Z}}\mbox{-}\underline{\rm
 mod} \; | \; {\rm
proj.dim}\;  {_AM_{n}} < \infty, \ \forall \ n\neq 0\}$$ and
$$K^b(A\mbox{-{\rm proj}})\simeq \mathcal{M}_P=\{\oplus_{n \in
\mathbb{Z}} M_{n}\in T(A)^{\mathbb{Z}}\mbox{-}\underline{\rm
 mod}\ | \; {\rm proj.dim}\;  {_AM_{n}} <
 \infty, \ \forall \ n \in \mathbb{Z}\}.$$
\end{thm}

\vskip10pt

\begin{cor}\ \ Let $A$ be a self-injective algebra. Then
$$D^b(A)\simeq\mathcal{N}=\{M=\oplus_{n \in \mathbb{Z}} M_{n}\in T(A)^{\mathbb{Z}}\mbox{-}\underline{\rm
 mod} \ | \;
{_AM_{n}} \ \mbox{is projective}, \ \forall \ n\neq 0\}$$ and
$$K^b(A\mbox{-{\rm proj}})\simeq\mathcal{M}_P=\{M=\oplus_{n \in
\mathbb{Z}} M_{n}\in T(A)^{\mathbb{Z}}\mbox{-}\underline{\rm
 mod}\ | \; {_AM_{n}} \mbox{ is projective}, \ \forall
\ n \in \mathbb{Z}\}.$$
\end{cor}

\vskip10pt

\subsection{}  For each $n\in \mathbb{Z}$ and an
indecomposable projective $A$-module $P$, the $\mathbb{Z}$-graded
$T(A)$-module
\begin{align} {\rm proj}(P, n, n+1) = \oplus_{i \in \mathbb{Z}} M_{i} \ \
\mbox{with} \ \ M_i = \begin{cases} P, \ \ &i= n;\\
D(A)\otimes_A P, \ \ & i =n+1;\\ 0, \ \
&\mbox{otherwise}.\end{cases}\end{align} is an indecomposable
projective $\mathbb{Z}$-graded $T(A)$-module, and any indecomposable
projective $\mathbb{Z}$-graded $T(A)$-module is of this form; and
for an indecomposable injective $A$-module $I$, the
$\mathbb{Z}$-graded $T(A)$-module
\begin{align}{\rm inj}(I, n-1, n) = \oplus_{i \in \mathbb{Z}} M_{i}\ \
\mbox{with} \ \ M_i = \begin{cases} {\rm Hom}_A(D(A), I), \ \ &i=
n-1;\\ I, \ \ & i =n;\\ 0, \ \
&\mbox{otherwise}.\end{cases}\end{align} is an indecomposable
injective $\mathbb{Z}$-graded $T(A)$-module, and any indecomposable
injective $\mathbb{Z}$-graded $T(A)$-module is of this form. Note
that
$${\rm proj}(P, n, n+1) \simeq {\rm inj}(D(A)\otimes_A P, n, n+1)$$and $${\rm inj}(I, n-1, n)
\simeq {\rm proj}({\rm Hom}_A(D(A), I), n-1, n).$$ Any homogeneous
$\mathbb{Z}$-graded $T(A)$-module $M=M_n$ of degree $n$ has the
injective hull ${\rm inj}(I_A(M_n), n-1, n)$, and the projective
cover ${\rm proj}(P_A(M_n), n, n+1)$, where $I_A(M_n)$ and
$P_A(M_n)$ are respectively the injective hull and the projective
cover of $M_n$ as an $A$-module. See [Hap1], II. 4.1.

\vskip10pt

\begin{lem} \ \ Let $A$ be a Gorenstein algebra. Then the full
subcategories given by
$$\{\oplus_{n \in \mathbb{Z}} M_{n}\in T(A)^{\mathbb{Z}}\mbox{-}\underline{\rm
 mod} \ | \; {\rm
proj.dim}\;  {_AM_{n}} < \infty, \ \forall\ n\neq 0\}$$ and
$$\{\oplus_{n \in \mathbb{Z}} M_{n}\in T(A)^{\mathbb{Z}}\mbox{-}\underline{\rm
 mod}\ | \; {\rm proj.dim}\;  {_AM_{n}} <
 \infty, \ \forall \ n \in \mathbb{Z}\}$$
are triangulated subcategories of $T(A)^{\mathbb Z}\mbox{-}
\underline{\rm mod}$.
\end{lem}
\noindent {\bf Proof.}\quad Since $A$ is Gorenstein, it follows from
$(3.2)$ that the two subcategories above contain all the injective
modules in $T(A)^{\mathbb{Z}}\mbox{-mod}$. Given a short exact
sequence $0 \longrightarrow M \longrightarrow N \longrightarrow L
\longrightarrow 0$ in $T(A)^{\mathbb{Z}}\mbox{-mod}$, then for each
$n$ we have an exact sequence of $A$-modules $0 \longrightarrow
M_{n} \longrightarrow N_{n} \longrightarrow
 L_{n} \longrightarrow 0$.
Note that  if any two terms of the short exact sequence above have
finite projective dimensions, then the other one also has finite
projective dimension. Now the assertion follows from Lemma 1.3.
\hfill $\blacksquare$

\par \vskip 10pt

\subsection{} {\bf Proof of Theorem 3.1.} \ \
We only prove
\begin {align}\mathcal{N} =  \{\oplus_{n \in
\mathbb{Z}} M_{n}\in T(A)^{\mathbb{Z}}\mbox{-}\underline{\rm
 mod} \ | \; {\rm proj.dim}\; {_AM_{n}} < \infty, \
\forall\ n\neq 0\}.\end{align} The other equality can be similarly
proven. Since $A$ is Gorenstein, it follows that
$$\{\oplus_{n \in \mathbb{Z}} M_{n} \ | \; {\rm proj.dim}\;
{_AM_{n}} < \infty, \ \forall\ n\neq 0\} = \{\oplus_{n \in
\mathbb{Z}} M_{n} \ | \; {\rm inj.dim}\; {_AM_{n}} < \infty, \
\forall\ n\neq 0\}.$$ By Lemma 3.3 the right hand side in $(3.3)$ is
a triangulated subcategory of $T(A)^{\mathbb Z}\mbox{-}
\underline{\rm mod}$ containing all the $A$-modules, while by
definition $\mathcal{N}$ is the triangulated subcategory of
$T(A)^{\mathbb Z}\mbox{-} \underline{\rm mod}$ generated by $A\mbox
{-\rm mod}$. It follows that $\mathcal{N}\subseteq  \{M=\oplus_{n
\in \mathbb{Z}} M_{n} \ | \; {\rm proj.dim}\; {_AM_{n}} < \infty, \
\forall\ n\neq 0\}.$ \vskip 10pt

For the other inclusion, first, consider all the objects of the form
$M=\oplus_{i\ge 0} M_{i}$ in the right hand side of $(3.3)$. We
claim that such an $M$ lies in $\mathcal{N}$, by using induction on
$l(M): = max\{\ i\ | \; M_{i} \neq 0\}$. Assume that $\mathcal{N}$
already contains all such objects $M$ with $l(M) < n$, \ $n\ge 1$.
Now, we use induction on $m: = {\rm inj.dim} \; _AM_{n}$ to prove
that $M=\oplus_{i = 0}^n M_{i}\in \mathcal{N},$ where ${\rm
inj.dim}\; {_AM_{i}} < \infty, \ \forall\ i\neq 0$.

\vskip10pt

If $m = 0$, i.e., $M_n$  is injective as an $A$-module, then by the
exact sequences in $TA^\mathbb{Z}\mbox{-}{\rm mod}$
\begin{align}0\longrightarrow M_n \longrightarrow M \longrightarrow M/{M_n} \longrightarrow 0.\end{align}
and (see $(3.2)$)
$$0 \longrightarrow M_n \longrightarrow {\rm inj}(M_n, n-1, n)
\longrightarrow M_n[1] \longrightarrow 0,$$  we have by induction
$M/{M_n}, \ {M_n}[1]\in \mathcal{N}$, and hence ${M_n}\in
\mathcal{N}$. By $(3.4)$ we have the distinguished triangle in
$T(A)^{\mathbb Z}\mbox{-} \underline{\rm mod}$
\begin{align} M_n \longrightarrow M \longrightarrow M/{M_n} \longrightarrow M_n[1]\end{align}
with  $M_n, \ M/{M_n}\in \mathcal{N}$. It follows from $\mathcal{N}$
being a triangulated subcategory that $M\in \mathcal{N}$.

\vskip10pt

Assume that for $n, d\ge 1$, $\mathcal{N}$ already contains all
the objects
 $M=\oplus_{i= 0}^n M_{i}$ in the right hand side of $(3.3)$ with ${\rm inj.dim} \;
_AM_{n} < d$. We will prove that $\mathcal{N}$ also contains such an
object $M$ with ${\rm inj.dim}\; {_AM_{n}}= d$. Take an exact
sequence in $T(A)^\mathbb{Z}\mbox{-mod}$ (see $(3.2)$)
$$0 \longrightarrow M_n \longrightarrow {\rm inj}(I_A(M_n), n-1, n) \longrightarrow M_n[1] \longrightarrow 0.$$
Since the $n$-th component $I_A(M_n)/M_n$ of $M_n[1]$ has injective
dimension less than $d$, it follows from induction that $M_n[1] \in
\mathcal{N}$, and hence $M_n \in \mathcal{N}$. Also $M/M_n\in
\mathcal{N}$ since $l(M/M_n) < n$. Thus $M\in \mathcal{N}$ by
$(3.5)$. This proves the claim.

\vskip10pt

Dually, any object of the form $M =\oplus_{i\le 0} M_{i}$ in the
right hand side of $(3.3)$ lies in $\mathcal{N}$.

\vskip10pt

In general, for  $M=\oplus_{n \in \mathbb{Z}} M_{n}$ in the right
hand side of $(3.3)$, set $M_{\geq 0}: = \oplus_{n \geq 0} M_{n}$.
Then it is a submodule of $M$. By the argument above we have
$M_{\geq 0}\in \mathcal{N}$ and $M/M_{\geq 0}\in \mathcal{N}$. Since
the exact sequence in $T(A)^\mathbb{Z}\mbox{-mod}$
$
0 \longrightarrow M_{\geq 0} \longrightarrow M \longrightarrow
M/M_{\geq 0} \longrightarrow 0
 $
induces a distinguished triangle in
$T(A)^\mathbb{Z}\underline{\mbox{-mod}}$, and $\mathcal{N}$ is a
triangulated subcategory, it follows that $M\in \mathcal{N}$. This
completes the proof. \hfill $\blacksquare$

\vskip 20pt

\section*{\bf Appendix: \  Stable category associated with
self-orthogonal modules}

\renewcommand{\thesection}{A}
\setcounter{subsection}{0} \setcounter{thm}{0}

We include a description of the stable category
$\underline{\mathfrak{a}(T)}$ of the Frobenius category
$\mathfrak{a}(T): =  \mathcal {X}_T \cap {_T}\mathcal {X}$  of
$T$-Cohen-Macaulay modules (cf. 2.1), where $T$ is a self-orthogonal
$A$-module. Denote by $K^{ac}(T)$ the full subcategory of the
unbounded homotopy category $K(A)$ consisting of acyclic complexes
with components in ${\rm add} T$. It is a triangulated subcategory.
Then we have

\vskip10pt

\begin{thm} Let  $T$ be a
self-orthogonal module such that ${\rm add}^\sim T\subseteq
T^\perp$ and ${\rm add}_\sim T \subseteq {^\perp T}$. Then there
is a triangle-equivalence $ \underline{{\mathfrak{a}(T)}} \simeq
K^{ac}(T)$. \end{thm}

\vskip5pt

We provide a direct proof. It can be also deduced by using Theorem
3.11 in [Bel]. Together with Theorem 2.5, we get another description
of the singularity category of a Gorenstein algebra. For a similar
result on separated noetherian schemes see [Kr], Theorem 1.1(3).
\vskip10pt

\begin{cor} \ \ Let $A$ be a Gorenstein algebra and $T$ be a generalized
tilting module. Then we have triangle-equivalences $\mathcal{D}_I(A)
\overset\sim\longleftarrow K^{ac}(T) \overset\sim\longrightarrow
\mathcal{D}_P(A).$
\end{cor}

\vskip10pt

\subsection{} Let $X\in {_T\mathcal{X}}$ with exact
sequence $$\cdots \longrightarrow T^{-i}
\stackrel{d_T^{-i}}{\longrightarrow} T^{-(i-1)} \longrightarrow
\cdots \stackrel{d_T^{-1}}{\longrightarrow} T^0
\stackrel{d_T^{0}}\longrightarrow X \longrightarrow 0,$$ where each
$T^{-i}\in {\rm add}T$ and ${\rm Ker}d^{-i} \in T^\perp$, $i\geq 0$.
Let $Y \in A\mbox{-mod}$ with a complex $$\cdots \longrightarrow
T'^{-i} \stackrel{d_{T'}^{-i}}{\longrightarrow} T'^{-(i-1)}
\longrightarrow \cdots \stackrel{d_{T'}^{-1}}{\longrightarrow} T'^0
\stackrel{d_{T'}^{0}}\longrightarrow Y \longrightarrow 0,$$ where
each $T'^{-i}\in {\rm add}T$. Denote them by $T^\bullet
\stackrel{d^0_T} \longrightarrow X$ and $T'^\bullet
\stackrel{d^0_{T'}} \longrightarrow Y$, respectively.

\vskip10pt

The proof of the following fact is similar with the one of the
Comparison-Theorem in homological algebra.

\vskip10pt

\begin{lem} With the notation of $X, Y, T^\bullet, T'^\bullet$ as above, for any morphism
$f: Y \longrightarrow X$, there exists a unique morphism $f^\bullet:
\ T'^\bullet \longrightarrow T^\bullet$ in $K(A)$ such that $f
d_{T'}^0= d_{T}^0 f^0$.
\end{lem}

\vskip10pt

The following fact is in p.446 in [Ric1], or p.45 in [KZ].

\vskip10pt

\begin{lem} Let $F: \mathcal {C}\longrightarrow \mathcal
{D}$ be a full and exact functor of triangulated categories. Then
$F$ is faithful if and only if it is faithful on objects, that is,
if $F(X)\simeq 0$ then $X\simeq 0$.
\end{lem} \par

\vskip10pt

\subsection{Proof of Theorem A.1.}\ \ Since
${\rm add}^\sim T\subseteq T^\perp$ and ${\rm add}_\sim T
\subseteq {^\perp T}$, it follows that $\ _T\mathcal {X} = {\rm
add}^\sim T$ and $\mathcal {X}_T = {\rm add}_\sim T.$ Thus, for
any object $T^\bullet$ in $K^{ac}(T)$ we have ${\rm Coker}
d_T^{i}\in {\rm add}^\sim T \cap {\rm add}_\sim T =
\mathfrak{a}(T)$, for each $i\in\Bbb Z$.

\vskip10pt Define a functor $F: K^{ac}(T) \longrightarrow
\underline{\mathfrak{a}(T)}$ as follows: for an object $T^\bullet$
in $K^{ac}(T)$, define $F(T^\bullet):={\rm Coker} d_T^{-1}$; for a
morphism $f^\bullet: T^\bullet \longrightarrow T'^\bullet$ in
$K^{ac}(T)$, define $F(f^\bullet)$ to be the image in
$\underline{\mathfrak{a}(T)}$ of the unique morphism $\bar{f^0}:
{\rm Coker}d_T^{-1} \longrightarrow {\rm Coker}d_{T'}^{-1}$ induced
by $f^0$. Note that $F$ is well-defined, dense, and full by  Lemma
A.3 and its dual.

\vskip10pt

Note that $F$ is faithful on objects. In fact, if
$F(T^\bullet)\simeq 0$, then ${\rm Coker}d_T^{-1} \in {\rm add}T$.
Since ${\rm Coker}d_T^{i} \in {\rm add}^\sim T\cap {\rm add}_\sim
T\subseteq T^\perp\cap \ ^\perp T$ \ for each $i$, it follows that
the exact sequence $0 \longrightarrow {\rm Coker} d_T^{-2}
\longrightarrow T^0 \longrightarrow {\rm Coker} d_T^{-1}
\longrightarrow 0$ splits, and hence  ${\rm Coker} d_T^{-2} \in
{\rm add}T$. Repeating this process, we have the split exact
sequence $0 \longrightarrow {\rm Coker} d_T^{-i-2} \longrightarrow
T^{-i} \longrightarrow {\rm Coker} d_T^{-i-1} \longrightarrow 0$,
and ${\rm Coker} d_T^{-i-2}\in {\rm add}T$, for each $i\ge 0$.
Similarly,  the exact sequence $0 \longrightarrow {\rm
Coker}d_T^{i-2} \longrightarrow T^i \longrightarrow {\rm
Coker}d_T^{i-1} \longrightarrow 0$ splits and ${\rm Coker}
d_T^{i-1} \in {\rm add}T$ for each $i\geq 1$. This implies that
the identity ${\rm Id}_{T^\bullet}$ is homotopic to zero, that is,
 $T^\bullet$ is zero in $K^{ac}(T)$.

\vskip10pt

In order to prove that $F$ is an exact functor, we first need to
establish a natural isomorphism $F\circ [1] \longrightarrow [1]\circ
F$, where the first $[1]$ is the usual shift of complexes, and the
second $[1]$ is the shift functor of the stable category
$\underline{\mathfrak{a}(T)}$. In fact, for each $T^\bullet \in
K^{ac}(T)$, we have a commutative diagram of short exact sequences
in $A\mbox{-mod}$

\[\xymatrix{
0 \ar[r] & F(T^\bullet) \ar[r]^{i_{T^\bullet}} \ar@{=}[d] & T^1
\ar[r]^{\pi_{T^\bullet}} \ar[d]^-{\gamma_{T^\bullet}} &
F(T^\bullet[1])
\ar[d]_{\alpha_{T^\bullet}} \ar[r]& 0\\
0 \ar[r] & F(T^\bullet) \ar[r]^-{i_{F(T^\bullet)}} & T(F(T^\bullet))
\ar[r]^{\pi_{F(T^\bullet)}} & F(T^\bullet)[1]\ar[r] & 0
 }\]
where  $i_{T^\bullet}$ is the natural embedding, $\pi_{T^\bullet}$
is the canonical map, and the second row is the one defining
$F(T^\bullet)[1]$, with $T(F(T^\bullet))\in {\rm add} T$ (see
$(1.1)$). Note that $\underline{\alpha_{T^\bullet}}$ is unique in
the stable category $\underline{\mathfrak{a}(T)}$, and that it is
easy to verify that $\underline{\alpha}: F\circ [1] \longrightarrow
[1]\circ F$ is a natural isomorphism (by using the same argument as
in the proof of Lemma 2.2 in [Hap1], p.12). We will show that $(F,
\underline{\alpha})$ is an exact functor.\par

\vskip10pt

Consider a distinguished triangle in $K^{ac}(T)$ by mapping cone
\begin{align*}
T^\bullet \stackrel{f^\bullet} \longrightarrow T'^\bullet
\stackrel{\binom{0}{1}}{\longrightarrow} {\rm
Con}(f^\bullet)\stackrel{(1, \; 0)} \longrightarrow T^\bullet [1].
\end{align*}
Write $\theta: = F(\binom {0}{1})$ and $\eta: = F((1, 0))$. Clearly
we have $\eta\theta = 0$. Observe that the sequence in $A$-mod
\begin{align*}
0 \longrightarrow F(T^\bullet)
\stackrel{\binom{\bar{f^0}}{-i_{T^\bullet}}}\longrightarrow
F(T'^\bullet)\oplus T^1 \stackrel{(\theta, \; \pi)} \longrightarrow
F({\rm Con}(f^\bullet)) \longrightarrow 0,
\end{align*}
is exact, where $\pi$ is the natural map from
$T^1$ to $F({\rm Con}(f^\bullet)) = (T^1\oplus T'^0)/ {\rm Im} \begin{pmatrix} -d_T^{0} & 0 \\
                            -f^{0} &  d_{T'}^{-1}      \end{pmatrix}$.
This can be seen by directly verifying that $(\theta, \pi)$ is
surjective, $(\theta, \pi)\binom{\bar{f^0}}{-i_{T^\bullet}}=0$, and
${\rm Ker}(\theta, \pi)\subseteq {\rm
Im}\binom{\bar{f^0}}{-i_{T^\bullet}}$. By definition we have $\eta
\pi=\pi_{T^\bullet}$, and hence the following diagram of short exact
sequences in $A$-mod commutes

\[\xymatrix{ 0 \ar[r] & F(T^\bullet)\ar@{=}[d]
\ar[r]^-{\binom{\bar{f^0}}{-i_{T^\bullet}}} & F(T'^\bullet)\oplus
T^1 \ar[d]^-{(0, -\gamma_{T^\bullet})} \ar[r]^-{(\theta, \; \pi)} &
F({\rm
Con}(f^\bullet)) \ar[d]^-{-(\alpha_{T^\bullet} \eta)} \ar[r] & 0\\
0 \ar[r] & F(T^\bullet) \ar[r]^-{i_{F(T^\bullet)}}& T(F(T^\bullet))
\ar[r]^-{\pi_{F(T^\bullet)}} & F(T^\bullet)[1] \ar[r] & 0.
 }\]

\vskip5pt  \noindent It follows from Lemma 1.2 that $F(T^\bullet)
\stackrel{F(f^\bullet)}\longrightarrow F(T'^\bullet)
\stackrel{\underline {\theta}}\longrightarrow F({\rm
Con}(f^\bullet)) \stackrel{ \underline{\alpha_{T^\bullet}}
\underline{\eta}}\longrightarrow F(T^\bullet)[1]$ is a distinguished
triangle in $\underline{\mathfrak{a}(T)}$. This proves that $F:
K^{ac}(T) \longrightarrow \underline{\mathfrak{a}(T)}$ is an exact
functor.

By Lemma A.4 the proof is completed. \hfill $\blacksquare$

\vskip10pt

\noindent{\bf Acknowledgement}: We are very grateful to the
anonymous referee for his or her helpful remarks and suggestions.

\vskip15pt

\bibliography{}

\end{document}